\documentclass[12pt]{article}%
\usepackage[applemac]{inputenc}
\usepackage{amsmath,amssymb}
\usepackage{amsfonts}
\usepackage{amsmath,amssymb}
\usepackage[applemac]{inputenc}
\usepackage{color}
\usepackage{color, colortbl}
\usepackage{amsmath}
\usepackage{amssymb}
\usepackage{bbm}
\usepackage{graphicx}
\usepackage{tikz}
\usepackage{geometry}
\usetikzlibrary{arrows}
\usepackage{amsfonts}
\usepackage{amsmath,amssymb}
\usepackage[applemac]{inputenc}
\usepackage{color}
\usepackage{amsmath}
\usepackage{amssymb}
\usepackage{graphicx}
\usepackage{multirow}
\usepackage{float}
\usepackage{caption}
\usepackage{subcaption}
\usepackage{bbm, dsfont}
\usepackage{tikz}
\newtheorem{theorem}{Theorem}[section]
\newtheorem{lemma}[theorem]{Lemma}
\newcommand{\E}{E}

\newtheorem{definition}[theorem]{Definition}

\newtheorem{example}[theorem]{Example}

\newtheorem{remark}[theorem]{Remark}

\usepackage{tikz}
\usetikzlibrary{shapes.geometric, arrows}

\tikzstyle{input} = [circle, minimum width=1cm, text centered, draw=black, fill=green!20]

\usetikzlibrary{arrows,shapes,automata,petri}

\newcommand{\dproof}{\noindent {Proof.} \quad}
\newcommand{\fproof}{\hfill $\square$ \bigskip}

\numberwithin{equation}{section}

\definecolor{LightCyan}{rgb}{0.88,1,1}

\def\1B{\text{1\!\!I}}

\def\<{\langle}
\def\>{\rangle}

\def\E{\mathbb{E}}
\def\R{\mathbb{R}}

\begin{document}
\title{Fokker-Planck equation for McKean-Vlasov SPDEs driven by time-space Brownian sheet}
\author{Nacira Agram$^{1},$ Bernt \O ksendal$^{2}$, Frank Proske$^{2}$ and Olena Tymoshenko$^{2, 3}$}
\date{29 April 2024
\vskip 0.8cm
Dedicated to the memory of Yuri Kondratiev}
\maketitle

\footnotetext[1]{Department of Mathematics, KTH Royal Institute of Technology 100 44, Stockholm, Sweden. \newline
Email: nacira@kth.se, pucci@kth.se.\\ Work supported by the Swedish Research Council grant (2020-04697) and the Slovenian Research and
Innovation Agency, research core funding No.P1-0448.}

\footnotetext[2]{%
Department of Mathematics, University of Oslo, Norway. \\
Emails: oksendal@math.uio.no, proske@math.uio.no, otymoshenkokpi@gmail.com}
\footnotetext[3]
{Department of Mathematical Analysis and Probability Theory, NTUU Igor Sikorsky Kyiv Polytechnic Institute, Kyiv, Ukraine.}

\begin{abstract}
In this paper, we consider a McKean-Vlasov (mean-field) stochastic partial differential equations (SPDEs) driven by a Brownian sheet. We study the propagation of chaos for a space-time Ornstein-Uhlenbeck SPDE type. Subsequently, we prove the existence and uniqueness of a nonlinear McKean-Vlasov SPDE. Finally, we establish a Fokker-Planck equation for the law of the solution of the McKean-Vlasov type SPDE driven by a time-space Brownian sheet, and we provide some examples to illustrate the results obtained.
\end{abstract}

%\tableofcontents

%BEGINNING
%\selectlanguage{english}
%%%%%%%%%%%%%%%%%%%%%%%%%%%%%%%%%%%%%%%%%%%%%%%%%%%%%%%%%%%%%%%%%%%%%%%%%%%%%%%%%%%%%%%%%%%%%%%%%%%%%%%%%%%%%%%%%%%%%%%%

\textbf{Keywords:} McKean-Vlasov (mean-field) SPDE; time-space Brownian sheet; Fokker-Planck equation.  
%%%%%%%%%%%%%%%%%%%%%%%%%%%%%%%%%
\section{Introduction}
In the one-parameter case, McKean-Vlasov equations were originally used in the kinetic theory of gases to model the dynamics of large particle systems in a medium. A crucial concept in kinetic theory is the propagation of chaos, which provides a description of the distribution of particles in a (monoatomic) gas medium as the number of particles tends to infinity. This concept, known as the "Stosszahlansatz," was assumed by Boltzmann (see \cite{Ehr}) to assert that each pair of particles in such a system moves independently. In this context, we also mention the work of Vlasov \cite{Vl}, which studied the propagation of chaos of charged particles in an electron gas or plasma. Regarding the investigation of the connection between Markov processes and nonlinear parabolic equations such as the Boltzmann equation in kinetic theory, we refer to McKean \cite{McK}.

%%%%%%%%%%%%%%%%%%%%%%%%%%%%%%%%%%%%%%%%%%%%%%%%%%%%%%%%%%%%%%%%%%%%%%%%%%%%%%%%%%%%%%%%%%%%%%%%%%%%%%%%%%%%%%%%%%%%%%%%%%%%%%%%%%%%%%%%%%%%%%%%%%%%%%%%%%%%%%%%%%%%%
The purpose of this paper is to study the law of the solution of McKean-Vlason type SPDEs driven by the Brownian sheet. We shall study the law $\mu_{t,x}=P_{Y(t,x)}$  of solutions $Y(t,x),t\geq
0,x\in\mathbb{R}$ of SPDEs of the form
\small
\begin{equation}
Y(t,x)=Y(t_{0},x_{0})+\int_{R(t,x)}\alpha(s,a,Y(s,a), \mu_{s,a})dsda
+\int_{R(t,x)}\beta(s,a,Y(s,a),\mu_{s,a})B(ds,da),\label{(1)}%
\end{equation}
where
$
R(t,x)=R^{(t_0,x_0)}(t,x)=[t_{0},t] \times[x_{0},x],t\geq t_{0},x\geq x_{0},
$
and $\alpha$ and $\beta$ are  Lipschitz continuous vector fields of linear growth,
$ B$ is a Brownian sheet.

The differential form of (\ref{(1)}) in terms of  time-space white noise  $ \overset{\bullet}{B}$ and Wick product $\diamond $ is
\small
\begin{equation}
\frac{\partial^{2}}{\partial t\partial x}Y(t,x)=\alpha(t, x,  Y
(t,x),\mu_{t,x})+\beta(t, x, Y(t,x),\mu_{t,x})\diamond\overset{\bullet}{B}(t,x).\label{(2)}
\end{equation}
The identity of (\ref{(1)}) and (\ref{(2)}) comes from the fact that
\small
\begin{equation*}
\int_{R(t,x)}\varphi(s,a)B(ds,da)=\int_{R(t,x)}\varphi(s,a
)\diamond\overset{\bullet}{B}(s,a)dsda,\text{ } \quad \text{for all} \quad \varphi, t, x.\label{(4)}%
\end{equation*}
See, for example, Holden et al. \cite{HOUZ} for more details.

The one-parameter case of McKean-Vlasov stochastic differential equations (SDEs) in infinite dimensions have been studied recently by Hong et al \cite{HLL}.

It is important to find the Fokker-Planck equation, denoted by $\mu_{t,x}$, for the McKean-Vlasov SPDE driven by a Brownian sheet as described above. It is worth noting that the specific type of Fokker-Planck equation for SPDEs driven by a Brownian sheet has not been previously addressed in the existing literature. Moreover, it is interesting to see the application of these concepts to understand the behavior of complex systems described by such equations.

However, in the one-parameter case, Fokker-Planck equations in infinite dimensions have been studied by Bogachev et al. \cite{BPR}. The authors developed a general technique to prove uniqueness of solutions for Fokker-Planck equations on infinite-dimensional spaces. They have been also studied by Agram et al. \cite{AO, APuO, AR} including various applications to optimal control and even deep learning.

Recently, Agram et al. \cite{AOPT} have considered various applications including the optimal control of time-space SPDEs driven by a Brownian sheet.

The paper is organized as follows: In Section 2, we review some preliminary concepts that will be used throughout this work. Specifically, we introduce background information about the stochastic calculus of time-space white noise. In Section 3, we investigate the propagation of chaos for a space-time Ornstein-Uhlenbeck SDE as a motivating factor for considering mean-field SDEs in the context of time-space. Section 4 is devoted to proving the existence and uniqueness of solutions to a McKean-Vlasov SPDE driven by a Brownian sheet. Finally, in Section 5, we state and prove the Fokker-Planck equation and illustrate the results for some time-space SPDE.

\section{Background}
In this section, we provide some background on the associated stochastic calculus for stochastic processes with two parameters. \\
Throughout this work, we denote by $\{B(t,x): t \geq 0, x \in \mathbb{R}\}$ a  Brownian sheet and $(\Omega, \mathcal{F}, P)$ a complete probability space on which we define the (completed) $\sigma$-field $\mathcal{F}_{t,x}$ generated by $ B(s,a), s \leq t, a \leq x$.
\begin{figure}[H]
\centering
\hspace{-0.5em}
\includegraphics[width=0.5\linewidth]{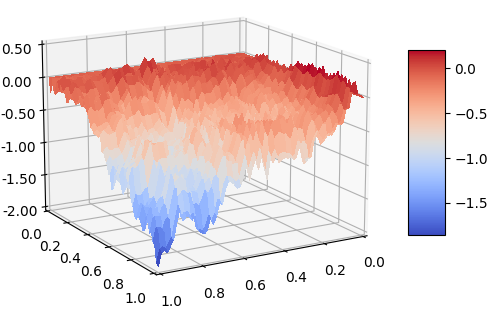}   
\caption{Brownian sheet on a 2D grid }
\end{figure}
Wong \& Zakai \cite{WZ} generalized the notion of stochastic integrals with respect to one-parameter Brownian motion to stochastic intergrals driven by the two-parameter Brownian sheet. Let us denote by $\mathbb{R}^2_+$ the positive quadrant of the plane, and let $z\in \mathbb{R}^2_+$.
 We define a first type stochastic integral with respect to the two-parameter Brownian sheet as introduced by Cairoli \cite{Cairoli72} denoted by:
\begin{align*}
   \int_{R_z} \phi(\zeta)B(d\zeta),
\end{align*}
and a second type \cite{WZ74} double stochastic integral denoted by
\begin{align*}
   \iint\limits_{R_{z}\times R_{z}}\psi(z,z')B(d\zeta)B(d\zeta'),
\end{align*}
where $R_z=[0,t]\times[0,x], z=(t,x), \zeta=(\zeta_1,\zeta_2)$ and $\zeta'=(\zeta_1',\zeta'_2)$.\\
In Wong \& Zakai \cite{WZ}, an It\^{o} formula for stochastic integrals in the plane is given.\\
In the sequel, we shall put sometimes $z=(t,x), \zeta=(s,a)$ in the following:\\
Let $\mathcal{P}$ be the predictable $\sigma$-algebra of subsets of $\Omega \times R_{z_0}$ generated by the sets $(z,z']\times A$, where $A\in \mathcal{F}_z$, and we denote by $\mathcal{D}$ the $\sigma$-algebra of $\Omega \times R_{z_0} \times R_{z_0}$ generated by the sets $(z_1,z_1']\times(z_2,z_2'] \times A$, where $(z_1,z_1']\bar{\wedge}
 (z_2,z_2']$ and $A\in \mathcal{F}_{z_1\vee z_2}$.
\subsection{The It\^{o} formula}
We will recall a two-parameter version of the It\^{o}
formula. First we introduce some notation
 from Wang \& Zakai \cite{WZ}. We shall introduce some notations which will be used throughout this work. We put $ \zeta=(\zeta_1,\zeta_2)=(s,a)\in \mathbb{R} \times \mathbb{R}$ and 
$d\zeta=d\zeta_{1}d\zeta_{2}=ds da$;
$B(t,x)$ is a Brownian sheet, $t\geq0,x\in\mathbb{R}$;
$z=(z_{1},z_{2})=(t,x),R_{z}=[0,z_{1}]\times\lbrack0,z_{2}]$;
$\int_{R_{z}}\varphi(\zeta)B(d\zeta)$ 
\text{denotes the It\^{o}
integral with respect to} 
$B(\cdot)$ \text{over} $R_{z}$;
$\int_{R_{z}}\psi(\zeta)d\zeta$ is two-dimensional Lebesgue integral of $\psi$; 
if $a=(a_{1},a_{2}),b=(b_{1},b_{2})$, then
$a\vee b=(\max(a_{1},b_{1}),\max(a_{2},b_{2})).
$
Moreover, 
\begin{equation*}
  \begin{aligned}
   I((a_1,a_2)\bar{\wedge} (b_1,b_2))=
\begin{cases}
        1 \quad & \text{if } a_1 \leq b_1\quad  \text{and } \quad a_2 \geq b_2, \\ 
        0 \quad & \text{otherwise. } 
    \end{cases}
  \end{aligned}
\end{equation*}

\begin{theorem}
[It\^{o} formula, Wang \& Zakai \cite{WZ}] 
Suppose
\begin{equation}
Y(z)=Y_{0}+\int_{R_{z}}\alpha(\zeta)d\zeta+\int_{R_{z}}\beta(\zeta
)B(d\zeta) + \iint\limits_{R_{z}\times R_{z}} \psi (\zeta,\zeta^{\prime})B(d\zeta)B(d\zeta^{\prime}).\label{(16)}%
\end{equation}
\end{theorem}
Then, if $f:\mathbb{R\rightarrow R}$ is smooth, we have
\small
\begin{align}
&f(Y(z))  =f(Y_{0})+\int_{R_{z}}f^{\prime}(Y(\zeta))[\alpha(\zeta
)d\zeta+\beta(\zeta)B(d\zeta)]\nonumber+\tfrac{1}{2}\int_{R_{z}}f^{\prime\prime}(Y(\zeta))\beta^{2}(\zeta
)d\zeta\nonumber\\
& +\iint\limits_{R_{z}\times R_{z}} \Big\{ f^{\prime\prime}(Y(\zeta\vee\zeta^{\prime
}))u \tilde{u}+f^{\prime} (Y(\zeta\vee \zeta^{\prime})) \psi(\zeta,\zeta^{\prime}) \Big\}B(d\zeta)B(d\zeta^{\prime}) \nonumber\\
&+\iint\limits_{R_{z}\times R_{z}}\Big\{f^{\prime\prime}(Y(\zeta\vee\zeta
^{\prime}))\Big( u \alpha(\zeta) + \psi(\zeta,\zeta^{\prime})  \tilde{u} \Big)\nonumber +\tfrac{1}{2}f^{(3)}(Y(\zeta\vee\zeta^{\prime}))u^2 \tilde{u}\Big\}d\zeta B(d\zeta^{\prime})\nonumber\\
&+\iint\limits_{R_{z}\times R_{z}}\Big\{f^{\prime\prime}(Y(\zeta\vee\zeta
^{\prime}))\Big( \tilde{u}\alpha(\zeta^{\prime}) + \psi(\zeta,\zeta^{\prime}) u  \Big)\nonumber +\tfrac{1}{2}f^{(3)}(Y(\zeta\vee\zeta^{\prime}))u^2 \tilde{u}\Big\}d\zeta B(d\zeta^{\prime})\nonumber\\
& +\iint\limits_{R_{z}\times R_{z}}I(\zeta \bar{\wedge} \zeta^{\prime}) \Big\{f^{\prime\prime}(Y(\zeta\vee\zeta
^{\prime}))\Big(\alpha(\zeta^{\prime})\alpha(\zeta)+\tfrac{1}{2}\psi^2(\zeta,\zeta^{\prime})\Big) \label{ito}\\
&+f^{(3)}(Y(\zeta\vee\zeta^{\prime}))u \tilde{u}\psi(\zeta,\zeta^{\prime})+\tfrac{1}{2}f^{(3)}(Y(\zeta\vee\zeta^{\prime}))\left[\alpha(\zeta^{\prime}%
)\tilde{u}^2+\alpha(\zeta)u^2 \right]\nonumber\\
& + \tfrac{1}{4}f^{(4)}(Y(\zeta\vee\zeta^{\prime}))u^2 \tilde{u}^2\Big\}d\zeta d\zeta^{\prime},\nonumber
\end{align}
%\end{theorem}

where $\displaystyle u=\beta(\zeta^{\prime})+\int_{R_{z}}I(\zeta \bar{\wedge} \zeta^{\prime})\psi(\zeta,\zeta^{\prime})B(d\zeta),$ $ \displaystyle \tilde{u}=\beta(\zeta)+\int_{R_{z}}I(\zeta \bar{\wedge} \zeta^{\prime})\psi(\zeta,\zeta^{\prime})B(d\zeta^{\prime}).$

\begin{remark}
Except for a deleted factor $\frac{1}{4}$ in the beginning of the term
\eqref{ito} this formula agrees with Proposition 5.1 in Wang \& Zakai \cite{WZ}. In the case $\alpha=0$ it is in agreement with the formula given by
Imkeller \cite{I1}, p. 35
\end{remark}
It is proved in \cite{WZ} that the double $B(d\zeta)B(d\zeta')$-integrals, and the mixed $d\zeta B(d\zeta')$ and $B(d\zeta)d\zeta'$-integrals are all weak martingales and hence have expectation $0$. Therefore, by the It\^{o} formula above we get the following:
\begin{theorem}
(Dynkin formula)
\small
\begin{align*}
 E[f(Y(z))]&   =f(Y_{0})+E\Big[\int_{R_{z}}\Big\{\alpha(\zeta)f^{\prime}(Y(\zeta
))+\frac{1}{2}\beta^{2}(\zeta)f^{^{\prime\prime}}(Y(\zeta))\Big\}d\zeta\\
&+\iint\limits_{R_{z}\times R_{z}}I(\zeta \bar{\wedge} \zeta^{\prime}) \Big\{f^{\prime\prime}(Y(\zeta\vee\zeta
^{\prime}))\Big(\alpha(\zeta^{\prime})\alpha(\zeta)+\tfrac{1}{2}\psi^2(\zeta,\zeta^{\prime})\Big)+f^{(3)}(Y(\zeta\vee\zeta^{\prime}))u \tilde{u}\psi(\zeta,\zeta^{\prime}) \\
&+\tfrac{1}{2}f^{(3)}(Y(\zeta\vee\zeta^{\prime}))\left[\alpha(\zeta^{\prime}%
)\tilde{u}^2+\alpha(\zeta)u^2 \right] + \tfrac{1}{4}f^{(4)}(Y(\zeta\vee\zeta^{\prime}))u^2 \tilde{u}^2\Big\}d\zeta d\zeta^{\prime}\Big].
\end{align*}
\end{theorem}

\begin{lemma}[Integration by parts]\label{partsa}
Suppose that for k=1,2
\small
\begin{align*}
    Y_k(z)=Y_k(0)+\int_{R_{z}}\alpha_k(\zeta)d\zeta
    +\int_{R_{z}}\beta_k(\zeta)B(d\zeta)+\iint\limits_{R_{z}\times R_{z}} \psi_k (\zeta,\zeta^{\prime})B(d\zeta)B(d\zeta^{\prime}).
\end{align*}
Then 
\small
\begin{align*}
  &  E[Y_1(z)Y_2(z)]=Y_1(0)Y_2(0)+E\Big[\int_{R_{z}}\Big\{ Y_1(\zeta) \alpha_2(\zeta)+Y_2(\zeta) \alpha_1(\zeta)+ \beta_1(\zeta)\beta_2(\zeta)\Big\}d\zeta\nonumber\\
& +\iint\limits_{R_{z}\times R_{z}}I(\zeta \bar{\wedge} \zeta^{\prime}) \Big\{\alpha_1 (\zeta^{\prime}) \alpha_2(\zeta) + \alpha_{1}(\zeta)\alpha_2(\zeta^{\prime}) +\psi_1(\zeta,\zeta^{\prime}) \psi_2(\zeta,\zeta^{\prime}) \Big\} d\zeta d\zeta^{\prime}\Big].\\
\end{align*}
\end{lemma}
%%%%%%%%%%%%%%%%%%%%%%%%%%%%%%%%%%%%%%%%%%%%%%5%%%%%%%%

\section{Propagation of Chaos for space time Ornstein-Uhlenbeck SDE}

In order to motivate mean-field SDEs in the two-parameter case from the
viewpoint of propagation of chaos, consider now the following linear $N-$%
particle system, $i=1,...,N$,
\begin{equation*}
Y^{i,N}(t,x)=Y_{i}(t_0,x_0)+\int_{0}^{t}\int_{0}^{x}%
\left(\frac{1}{N}\sum_{j=1}^{N}a_{j}Y^{j,N}(\zeta_{1},\zeta_{2})-Y^{i,N}(\zeta_{1},\zeta_{2})\right)d\zeta_{1}d\zeta_{2}+B^{i}(t,x),
\end{equation*}%
where  $a_{j}\in \mathbb{R}$, $%
j=1,...,N$ and $B=(B^{1},...,B^{N})$ a $N-$dimensional Brownian sheet.

So the latter system of equations, which could be e.g. used to describe the
dynamics of interacting waves in an ocean, can be written as
\begin{equation*}
Y(t,x)=Y(t_0,x_0)+\int_{0}^{t}\int_{0}^{x}\left(\frac{1}{N}%
A-I_{N}\right)Y(\zeta_{1},\zeta_{2})d\zeta_{1}d\zeta_{2}+B(t,x),
\end{equation*}%
where  $0\leq t,x\leq T$, $Y(t_0,x_0)\in \mathbb{R}%
^{N}$, $Y(t,x):=(Y^{1,N}(t,x),...,Y^{N,N}(t,x))$, $I_{N}$ is the unit matrix
and%
\begin{equation*}
A:=
\begin{pmatrix} a_{1} & \cdots  & a_{N} \\ 
\vdots  &  & \vdots  \\ 
a_{1} & \cdots  & a_{N} \end{pmatrix}.
\end{equation*}
By applying the Malliavin derivative $D_{u,v}$ to both sides of the latter
equation, for $0\leq
u\leq t$, $0\leq v\leq x$, we obtain that%
\begin{equation*}
D_{u,v}Y(t,x)=\int_{u}^{t}\int_{v}^{x}\left(\frac{1}{N}%
A-I_{N}\right)D_{u,v}Y(\zeta_{1},\zeta_{2})d\zeta_{1}d\zeta_{2}+\chi _{R(t,x)}(u,v)I_{N}.
\end{equation*}%
Using Picard iteration we see that%
\begin{equation*}
   D_{u,v}Y(t,x) =\sum_{n\geq
0}\int_{R(u,v,t,x)}\int_{R(u,v,t_{1},x_{1})}...\int_{R(u,v,t_{n-1},x_{n-1})}
\prod\limits_{j=1}^{n}\left(\frac{1}{N}A-I_{N}\right)dt_{n}dx_{n}...dt_{1}dx_{1} 
\end{equation*}
\begin{equation*}
    =\sum_{n\geq 0}\frac{1}{(n!)^{2}}\left( (t-u)(x-v)t\left(\frac{1}{N}A-I_{N}\right)\right)^{n}=f\left((t-u)(x-v)\left(\frac{1}{N}A-I_{N}\right)\right)\text{,}
\end{equation*}
where $f$ is a function, which is related to the Bessel function of order
zero and given by  
$\displaystyle
f(y)=\sum_{n\geq 0}\frac{1}{(n!)^{2}}y^{n}\text{.}
$
Then it follows from the Clark-Ocone theorem for a Brownian sheet that%
\begin{equation*}
Y(t,x)=E\left[ Y(t,x)\right] +\int_{0}^{t}\int_{0}^{x}f\left((t-u)(x-v)\left(\frac{1}{N}A-I_{N}\right)\right)B(du,dv)\text{.}
\end{equation*}%
Similarly, we find that%
\begin{equation*}
E\left[ Y(t,x)\right] =f\left(tx\left(\frac{1}{N}A-I_{N}\right)\right)Y(t_0,x_0)\text{.}
\end{equation*}%
Hence,%
\begin{equation*}
Y(t,x)=f\left(tx\left(\frac{1}{N}A-I_{N}\right)\right)Y(t_0,x_0)+\int_{0}^{t}\int_{0}^{x}f\left((t-u)(x-v)\left(\frac{1}{N}A-I_{N}\right)\right)B(du,dv)\text{.}
\end{equation*}%
Define 
$
\displaystyle \left\Vert A\right\Vert =\sum_{j=1}^{N}a_{j}\text{.}
$ Require that $Y_{i}(t_0,x_0)=y$ for all $i=1,...,N$, the sequence $a_{j},j\geq 1$ is
bounded and that%
\begin{equation*}
\frac{1}{N}\left\Vert A\right\Vert =\frac{1}{N}\sum_{j=1}^{N}a_{j}\underset{%
N\longrightarrow \infty }{\longrightarrow }a>0\text{.}
\end{equation*}%
On the other hand, we observe for $n\geq 0$ that%
\begin{eqnarray*}
\left(\frac{1}{N}A-I_{N}\right)^{n} &=&\sum_{j=0}^{n}\binom{n}{j}\left(\frac{1}{N}A\right)^{j}(-I_{N})^{n-j} =(-1)^{n}I_{N}+\sum_{j=1}^{n}\binom{n}{j}\left(\frac{1}{N}A\right)^{j}(-I_{N})^{n-j}
\\
&=&(-1)^{n}I_{N}+\sum_{j=1}^{n}\binom{n}{j}\left\Vert \frac{1}{N}%
A\right\Vert ^{j-1}\frac{1}{N}A(-I_{N})^{n-j} \\
&=&(-1)^{n}I_{N}+\left(-(-1)^{n}+(-1)^{n}+\sum_{j=1}^{n}\binom{n}{j}\left\Vert 
\frac{1}{N}A\right\Vert ^{j}(-1)^{n-j}\right)\frac{1}{\left\Vert \frac{1}{N}%
A\right\Vert }\frac{1}{N}A \\
%&=&(-1)^{n}\left(I_{N}+(-1+(-1)^{n}(\frac{1}{N}\left\Vert A\right\Vert -1)^{n}))%
%\frac{1}{\left\Vert A\right\Vert }A\right) \\
&=&(-1)^{n}(I_{N}-\frac{1}{\left\Vert A\right\Vert }A)+(\frac{1}{N}%
\left\Vert A\right\Vert -1)^{n}\frac{1}{\left\Vert A\right\Vert }A\text{.}
\end{eqnarray*}%
So%
\begin{eqnarray*}
&&f\left((t-u)(x-v)(\frac{1}{N}A-I_{N})\right)
=\sum_{n\geq 0}\frac{1}{(n!)^{2}}((t-u)(x-v)(\frac{1}{N}A-I_{N}))^{n} \\
&=&\sum_{n\geq 0}\frac{1}{(n!)^{2}}((t-u)(x-v))^{n}\left\{ (-1)^{n}(I_{N}-%
\frac{1}{\left\Vert A\right\Vert }A)+(\frac{1}{N}\left\Vert A\right\Vert
-1)^{n}\frac{1}{\left\Vert A\right\Vert }A\right\}  \\
&=&f(-(t-u)(x-v))(I_{N}-\frac{1}{\left\Vert A\right\Vert }A)+f((t-u)(x-v)(%
\frac{1}{N}\left\Vert A\right\Vert -1))\frac{1}{\left\Vert A\right\Vert }A%
\text{.}
\end{eqnarray*}%
Hence%
\begin{eqnarray*}
Y(t,x) &=&f(-tx)(I_{N}-\frac{1}{\left\Vert A\right\Vert }A)y+f(tx)\left(\frac{1}{N}%
\left\Vert A\right\Vert -1\right)\frac{1}{\left\Vert A\right\Vert }Ay \\
&&+\int_{0}^{t}\int_{0}^{x}\left\{ f(-(t-u)(x-v))\left(I_{N}-\frac{1}{\left\Vert
A\right\Vert }A\right)\right.  \\
&&\left. +f((t-u)(x-v)\left(\frac{1}{N}%
\left\Vert A\right\Vert -1\right)\frac{1}{%
\left\Vert A\right\Vert }A\right\} B(du,dv)\text{.}
\end{eqnarray*}%
The latter entails that%
\begin{eqnarray*}
X^{i,N}(t,x) &=&f(-tx)\left(1-\frac{\left\Vert A\right\Vert }{\left\Vert
A\right\Vert }\right)y+f(tx)\left(\frac{1}{N}%
\left\Vert A\right\Vert -1\right)y \\
&&+\int_{0}^{t}\int_{0}^{x}f(-(t-u)(x-v))dB^{i}(du,dv)+I_{i,N}\text{,}
\end{eqnarray*}%
where 
\begin{eqnarray*}
I_{i,N} &:=&-\sum_{j=1}^{N}\int_{0}^{t}\int_{0}^{x}f(-(t-u)(x-v))\frac{a_{j}%
}{\left\Vert A\right\Vert }dB^{j}(du,dv) \\
&+&\sum_{j=1}^{N}\int_{0}^{t}\int_{0}^{x}f(-(t-u)(x-v))\left(\frac{1}{N}%
\left\Vert A\right\Vert -1\right)\frac{a_{j}}{\left\Vert A\right\Vert }
dB^{j}(du,dv),
\end{eqnarray*}
for all $1\leq i\leq N$. Using the It\^{o} isometry, our assumptions on $%
a_{j},j\geq 1$ and dominated convergence, we find that%
\begin{eqnarray*}
E\left[ \left\vert I_{i,N}\right\vert ^{2}\right] 
&\leq &C\frac{N}{\left\Vert A\right\Vert ^{2}}\int_{0}^{t}%
\int_{0}^{x}(f(-(t-u)(x-v)))^{2}(\frac{1}{N}\left\Vert A\right\Vert
-2)^{2}dudv \\
&&\underset{N\longrightarrow \infty }{\longrightarrow }0\text{,}
\end{eqnarray*}
where $C$ is a constant.
Thus,  for $N\longrightarrow \infty$ 
\begin{equation*}
Y^{i,N}(t,x)\overset{d}{\longrightarrow }Y(t,x):=f(tx)(a-1)y+\int_{0}^{t}%
\int_{0}^{x}f(-(t-u)(x-v))dB(du,dv),
\end{equation*}%
for all $i\geq 1$, where $Y(t,x)$, $0\leq t$, $x\leq T$ solves (by the same
reasoning as above) the mean-field hyperbolic SPDE  
\begin{equation*}
Y(t,x)=y+\int_{0}^{t}\int_{0}^{x}aE\left[ Y(\zeta_{1},\zeta_{2})\right]
-Y(\zeta_{1},\zeta_{2})d\zeta_{1}d\zeta_{2}+B(t,x).
\end{equation*}%

%%%%%%%%%%%%%%%%%%%%%%%%%%%%%%%%%%%%%%%%%%%%
\section{McKean-Vlasov SPDE}
In this section, we prove the existence and uniqueness of the solution to a McKean-Vlasov SPDE driven by a Brownian sheet. To this end, we denote by $z=(t,x)\in \mathbb{R}^{2}_+$, ${R}_z$ will denote the rectangle $[0,t]\times[0,x]$. The object of our study is an SPDE of the type
\begin{equation} \label{MV}
    \frac{\partial^2}{\partial t \partial x} Y(t,x)=\alpha(t,x,Y(t,x),\mu_{t,x})+\beta(t,x,Y(t,x),\mu_{t,x})\diamond \overset{\bullet}{B}(t,x),\quad
Y(0,0)=y.
\end{equation}
Here $\mu_{t,x}$ represents the probability distribution of $Y(t,x)$. We shall define the space where these probability distributions reside.
\begin{definition}[Special weighted Sobolev space]
Let $\mathbb{M}$ be the pre-Hilbert space of random measures $\mu$ on
$\mathbb{R}$ equipped with the norm  %\textcolor{red}{do we need expectation?}
\begin{align}\label{norm}
\left\Vert \mu\right\Vert _{\mathbb{M}}^{2}:= \mathbb{E[}\int_{\mathbb{R}}
|\hat{\mu}(y)|^{2}e^{-y^{2}}dy]\text{,}
\end{align}
where $\hat{\mu}$ is the Fourier transform of the measure $\mu$, i.e.%
\[%
\begin{array}
[c]{lll}%
\hat{\mu}(y) & := & {\int_{\mathbb{R}}
}e^{-ixy}\mu(dx);\quad y\in\mathbb{R}.
\end{array}
\]
If $\mu,\eta\in
\mathbb{M}$, we define the inner product $\left\langle \mu
,\eta\right\rangle _{\mathbb{M}}$ by
\small
\[
\left\langle \mu,\eta\right\rangle _{\mathbb{M}}=\E{[}
\int_{{\mathbb{R}}}\operatorname{Re}(\overline{\hat{\mu}}(y)\hat{\eta
}(y))e^{-y^{2}}dy]  ,
\] 
\normalsize
where, $\operatorname{Re}(z)$ denotes the real part and $\bar{z}$ denotes the complex conjugate of the complex number $z$. 
\end{definition}
The space $\mathbb{M}$ equipped with the inner
product $\left\langle \mu,\eta\right\rangle _{\mathbb{M}}$
is a pre-Hilbert space. For not having ambiguity, we will also use the notation $\mathbb{M}$ for the completion of this pre-Hilbert space.
%%%%%%%%%%%%%%%%%%%%%%%%%%%
Moreover, we have the following estimate:
Let $Y_{1}$ and $Y_{2}$ be two $d$-dimensional random variables in
$L^{2}(\mathbb{P})$ with associated probability distributions $\mu_{1}$ and $\mu_{2}$, respectively. Thus
\[
\begin{array}
[c]{lll}%
\left\Vert \mu_{1}-\mu_{2}\right\Vert _{\mathbb{M}
}^{2} & \leq & \pi\ \mathbb{E}[(Y_{1}-Y_{2})^{2}]\text{.}%
\end{array}
\]
To study the well-posedness of the McKean-Vlasov SPDE  \eqref{MV} driven by Brownian sheet, we impose the following set of assumptions on the coefficients $\alpha$ and $\beta$ for which we aim, they will insure the well-posedness:\\
\begin{itemize}
    \item [(a)] $\alpha(z,y,\mu):[0,T]^2\times \mathbb{R} \times \mathbb{M}\rightarrow \mathbb{R}, \beta(z,y,\mu):[0,T]^2\times \mathbb{R} \times \mathbb{M}\rightarrow \mathbb{R}$ are locally bounded and Borel-measurable functions.
      \item  [(b)] There exists a constant $C$, such that for all $z\in \mathbb{R}^{2}_+,y,y',\mu,\mu'$, we have
    \begin{align*}
       \left\vert \alpha(z,y,\mu)-\alpha(z,y',\mu')\right\vert &+\left\vert \beta(z,y,\mu)-\beta(z,y',\mu')\right\vert \leq C (\left\vert y-y'  \right\vert  +\left\vert \left\vert \mu-\mu'  \right\vert\right\vert _{\mathbb{M}}),
    \end{align*} 
and 
\begin{align*} 
\E\int_{{R}_z}\left\vert \alpha(\zeta,0,\delta_{(0,0)})\right\vert^2+\left\vert \beta(\zeta,0,\delta_{(0,0}))\right\vert^2 d\zeta\leq \infty,
\end{align*} 
where $\delta$ is the Dirac measure.
\end{itemize}
Let $f_0$ be the Bessel function of order zero and $r_0\approx 1.4458$ be the first nonnegative zero of $J_0$:
$$
r_0=inf \left\{t>0:f_0(2\sqrt{t})=\sum_{j=0}^\infty \frac{(-1)^j}{{j!}^2} t^j=0 \right\}.
$$
We shall recall the two-parameter version of Gronwall's
Lemma in \cite{ZN}.

\begin{lemma}[Two-parameter Gronwall's Lemma]
   Let $f$ be a non-negative and bounded function. There exists $C_0>0$ satisfies $C_0|z|\leq r_0$, such that
   $$f(z) \leq C_0 \int_{{R}_z} f(\zeta) d\zeta. $$
   Moreover, $f$ vanishes on $R_z$.
\end{lemma}
\begin{theorem}[Existence and uniqueness]
Under the above assumptions (a)-(b), the conditional McKean-Vlasov SPDE \eqref{MV} has a unique strong solution. 
\end{theorem}
\dproof 
The proof is based on the Picard iteration argument as in the proof of the Propagation of Chaos for Space time Ornstein-Uhlenbeck SPDE in Section 3.\\
\textbf{Step 1. Uniqueness}\\
Suppose that we have two solutions $Y,Y'$ and set $\tilde{Y}=Y-Y'$, such that $\tilde{Y}$ satisfies
\begin{align*}
\tilde{Y}(t,x)&=\int_{{R}_z}\{\alpha(\zeta,Y(\zeta),\mu_{\zeta})-\alpha(\zeta,Y'(\zeta),\mu'_{\zeta})\}d\zeta\nonumber\\
&+\int_{{R}_z}\{\beta(\zeta,Y(\zeta),\mu_{\zeta})-\beta(\zeta,Y'(\zeta),\mu'_{\zeta})\}dB(\zeta).
\end{align*}
Taking the mean square yields
\begin{align*}
E|\tilde{Y}(t,x)|^2&=E|\int_{{R}_z}\{\alpha(\zeta,Y(\zeta),\mu_{\zeta})-\alpha(\zeta,Y'(\zeta),\mu'_{\zeta})\}d\zeta\nonumber\\
&+\int_{{R}_z}\{\beta(\zeta,Y(\zeta),\mu_{\zeta})-\beta(\zeta,Y'(\zeta),\mu'_{\zeta})\}dB(\zeta)|^2.
\end{align*}
Triangle inequality together with the linearity of the expectation, leads
\begin{align*}
E|\tilde{Y}(t,x)|^2&\leq E|\int_{{R}_z}\{\alpha(\zeta,Y(\zeta),\mu_{\zeta})-\alpha(\zeta,Y'(\zeta),\mu'_{\zeta})\}d\zeta|^2\nonumber\\
&+E|\int_{{R}_z}\{\beta(\zeta,Y(\zeta),\mu_{\zeta})-\beta(\zeta,Y'(\zeta),\mu'_{\zeta})\}dB(\zeta)|^2.
\end{align*}
We use the Cauchy-Schwarz inequality for the $d\zeta$-integral and the isometry for the $dB$-integral, we get 
\begin{align*}
E|\tilde{Y}(t,x)|^2&\leq |z|^2 E\int_{{R}_z}|\alpha(\zeta,Y(\zeta),\mu_{\zeta})-\alpha(\zeta,Y'(\zeta),\mu'_{\zeta})|^2d\zeta\nonumber\\
&+E\int_{{R}_z}|\beta(\zeta,Y(\zeta),\mu_{\zeta})-\beta(\zeta,Y'(\zeta),\mu'_{\zeta})|^2d\zeta,
\end{align*}
where we have used the notation $|z|=tx.$\\
Jensen inequality combined with the Lipschitz condition, give
\begin{align*}
E|\tilde{Y}(t,x)|^2&\leq |z|^2 (C+C\pi)^2\int_{{R}_z} E|Y(\zeta)-Y'(\zeta)|^2d\zeta\nonumber+(C+C\pi)^2\int_{{R}_z} E|Y(\zeta)-Y'(\zeta)|^2d\zeta\nonumber\\
&=(C+C\pi)^2(|z|^2+1)\int_{{R}_z} E|\tilde{Y}(\zeta)|^2d\zeta.
\end{align*}
By Gronwall's Lemma , we get
$E|\tilde{Y}(t,x)|^2=0.$\\
\textbf{Step 2. Existence}\\
Define $Y^0(z)=y$ and $Y^n(z)$ inductively with corresponding probability distributions $\delta_y$ and $\mu^{n}_{\zeta}=P_{Y^n(\zeta)}$ respectively, as follows
\begin{align*}
Y^{n+1}(z)&=y+\int_{{R}_z}\alpha(\zeta,Y^{n}(\zeta),\mu^{n}_{\zeta})d\zeta+\int_{{R}_z}\beta(\zeta,Y^{n}(\zeta),\mu^{n}_{\zeta})dB(\zeta).
\end{align*}
Similar computations as in the uniqueness, for some constant $K>0$ depending on the lipschitz constant, lead to
\begin{align*}
E|Y^{n+1}(z)-Y^n(z)|^2&\leq K^2|z|^2\int_{{R}_z} E|Y^{n}(\zeta)-Y^{n-1}(\zeta)|^2d\zeta.
\end{align*}
Repeating this procedure $n$-times, we get 
\small
\begin{align*}
&E|Y^{n+1}(z)-Y^n(z)|^2\\
&\leq K^{2n}|z|^2\int_{{R}_{z}}\int_{{R}_{z_{n}}} \cdot \cdot \cdot (\int_{{R}_{z_{n}}} E|Y^{1}(\zeta_{n+1})-Y^{0}(\zeta_{n+1})|^2d\zeta_{n+1})d\zeta_{{n}}
  \cdot \cdot \cdot d\zeta_{z_1}\\
 &\leq K^{2n}|z|^{2n} sup_{u \in {R}_{z}} E|Y^{1}(u)|^2 x_n.
  \end{align*}
  Taking the sum, we have
  \small
\begin{align*}
\sum_{n=0}^
{\infty} E|Y^{n+1}(z)-Y^n(z)|^2&\leq sup_{u \in {R}_{z}} E|Y^{1}(u)|^2 \sum_{n=0}^
{\infty} (K|z|)^{2n} x_n<\infty,
\end{align*}
with $x_n=-\sum_{j=1}^
{n}\frac{(-1)^j}{(j!)^{2}}x_{n-j}$ and $K|z|<\sqrt{r_0}$. Thus $(Y^n)_n$ converges.
\fproof

%%%%%%%%%%%%%%%%%%%%%%%%%%%%%%%%%%%%%%%%%%%%
\section{The Fokker-Planck equation}
In this section we state and prove our main result.\\
In the following we let $D=\frac{\partial}{\partial y}$ denote the derivative in the sense of distributions meaning in the space $\mathcal{S}'$ of tempered distributions on $\mathbb{R}$.\\

We first prove the following auxiliary result:
\begin{lemma} \label{4.1}
Suppose $F$ has the form
\small
    \begin{align}
       F(\zeta,\zeta',y,\mu)= D^p[f(\zeta,y,\mu)g(\zeta',y,\mu)\mu_{\zeta \vee \zeta'}];\quad p=1,2,3,4.
    \end{align}
Then 
\small
    \begin{align*}
        \frac{\partial^2}{\partial t \partial x} \int_{R_z} \int_{R_z} I(\zeta \bar{\wedge} \zeta')F(\zeta,\zeta',y,\mu)d\zeta d\zeta'=D^p\Big[\Big(\int_0^t  f((\zeta_1,x),y,\mu)d\zeta_1 \int_0^xg((t,\zeta'_2),y,\mu)) d\zeta'_2\Big)\mu_{t,x}\Big].
    \end{align*}
In particular, if $f(\zeta,y,\mu)=f(\zeta)$ and $g(\zeta,y,\mu)=g(\zeta)$do not depend on $y$ and $\mu$, we get
\small
  \begin{align*}
        \frac{\partial^2}{\partial t \partial x} \int_{R_z} \int_{R_z} I(\zeta \bar{\wedge} \zeta')F(\zeta,\zeta',y,\mu)d\zeta d\zeta'=\Big(\int_0^t  f(\zeta_1,x)d\zeta_1 \int_0^x g(t,\zeta'_2) d\zeta'_2\Big)D^p\mu_{t,x}, 
    \end{align*}  
and if $f,g$ are a constants, we have
\small
  \begin{align*}
        \frac{\partial^2}{\partial t \partial x} \int_{R_z}\int_{R_z} I(\zeta \bar{\wedge} \zeta')F(\zeta,\zeta',y,\mu)d\zeta d\zeta'=txfg D^p\mu_{t,x}. 
    \end{align*}
    \end{lemma}
\dproof
By the definition of $I(\zeta \bar{\wedge} \zeta')$ we have
\begin{align*}
    H:=\int_{R_z}\int_{R_z} I(\zeta \bar{\wedge} \zeta')F(\zeta,\zeta') d\zeta d\zeta'&= \int_{R_z}\int_0^x \int_0^t \mathbf{1}_{\zeta_1 \leq \zeta'_1 \& \zeta_2 \geq \zeta'_2} F(\zeta,\zeta') d\zeta d\zeta' \nonumber\\
    &= \int_0^x \int_0^t \left(\int_0^{\zeta_2} \int_{\zeta_1}^t F(\zeta,\zeta') d\zeta'_1 d\zeta'_2\right) d\zeta_1 d\zeta_2.
\end{align*}
Taking derivatives, we get
\begin{align*}
\frac{\partial }{\partial t }H= \int_0^x \int_0^t \left(\int_0^{\zeta_2} F(\zeta,(t,\zeta'_2))d\zeta'_2\right) d\zeta_1 d\zeta_2,
\end{align*}
and
\begin{align*}
\frac{\partial^2 }{\partial x \partial t }H&= \int_0^t \int_0^x  F((\zeta_1,x),(t,\zeta'_2))d\zeta'_2) d\zeta_1\nonumber= \int_0^t \int_0^x D^p[f(\zeta_1,x) g(t,\zeta'_2)\mu_{(\zeta',x) \vee (t,\zeta'_2)}]\nonumber\\
&= \int_0^t \int_0^x D^p[f(\zeta_1,x) g(t,\zeta'_2)\mu_{t,x}]d\zeta'_2 d\zeta_1\nonumber= D^p\Big[\Big(\int_0^t \int_0^x f(\zeta_1,x) g(t,\zeta'_2)d\zeta'_2 d\zeta_1\Big)\mu_{t,x}\Big]\nonumber\\
&= D^p\Big[\Big(\int_0^t  f(\zeta_1,x) d\zeta_1 \int_0^x g(t,\zeta'_2)d\zeta'_2\Big)\mu_{t,x}\Big]
\end{align*}
\fproof

\begin{theorem}[The Fokker-Planck equation]
Let $Y(t,x)$ be the solution of a mean-field SPDE of the form
\begin{align*}
\frac{\partial^2}{\partial t \partial x} Y(t,x)&=\alpha(t,x,Y(t,x),\mu_{t,x})+\beta(t,x,Y(t,x),\mu_{t,x})\diamond \overset{\bullet}{B}(t,x), \,\, t, x>0;\,\, 
Y(0,0)=y.
\end{align*}
\begin{itemize}
\item
(Integral version) Then the law $\mu_{t,x}=\mu_{t,x}(dy)=\mathcal{L}(Y(t,x))(dy)$ of $Y(t,x)$ satisfies the following integral  (Fokker-Planck)  equation:

\begin{align}
    \mu_z = \delta_y + \int_{R_z}\Big\{A_1[\mu_{\zeta}] + \int_{R_z}L_{\zeta'}[\mu_{\zeta}] d\zeta' \Big \} d\zeta,\label{FP1}
\end{align}
where the operators $A_1$ and $L_{\zeta'}$ are given by
\begin{align}\label{A1}
A_1[\mu_{\zeta}] = -D[\alpha(\zeta,y,\mu) \mu_{\zeta}] + \tfrac{1}{2} D^2[\beta^2(\zeta,y,\mu_{\zeta}) \mu_{\zeta}]
\end{align}
and, with $z=(t,x)$,
\begin{align}\label{L}
L_{\zeta'}[\mu_{\zeta}] &= I(\zeta \bar{\wedge} \zeta^{\prime}) \Big( D^2\Big[\alpha(\zeta,y,\mu)\alpha(\zeta',y,\mu) \mu_{\zeta \vee \zeta'}\Big] \nonumber\\
&- \tfrac{1}{2} D^3\Big[\Big(\alpha(\zeta,y,\mu)\beta^2(\zeta',y,\mu)+\alpha(\zeta',y,\mu)\beta^2(\zeta,y,\mu)\Big) \mu_{\zeta \vee \zeta'}\Big]\nonumber\\
& +\tfrac{1}{4} D^4\Big[ \beta^2(\zeta,y,\mu) \beta^2(\zeta',y,\mu)\mu_{\zeta \vee \zeta'} \Big]\Big).
\end{align}
\item
(Differential form) Equivalently, in differential form the Fokker-Planck equation states that
\begin{align}
    \frac{\partial^2}{\partial t \partial x} \mu_{t,x}&= -D[\alpha(z,y,\mu) \mu_{z}] + \tfrac{1}{2} D^2[\beta^2(z,y,\mu) \mu_{z}]\nonumber\\
    &+ D^2\Big[\Big(\int_0^t \alpha((\zeta_1,x),y,\mu)d\zeta_1 \int_0^x\alpha((t,\zeta'_2),y,\mu)d\zeta'_2\Big) \mu_{t,x}\Big] \nonumber\\
&- \tfrac{1}{2} D^3\Big[\Big(\int_0^t \alpha((\zeta_1,x),y,\mu)d\zeta_1 \int_0^x\beta^2((t,\zeta'_2),y,\mu)d\zeta'_2 \nonumber\\
&+\int_0^x\alpha((t,\zeta'_2),y,\mu)d\zeta'_2 \int_0^t \beta^2((\zeta_1,x),y,\mu)d\zeta_1\Big) \mu_{t,x}\Big]\nonumber\\
& +\tfrac{1}{4} D^4\Big[\Big(\int_0^t  \beta^2((\zeta_1,x),y,\mu)d\zeta_1 \int_0^x \beta^2((t,\zeta'_2),y,\mu) d\zeta'_2\Big)\mu_{t,x} \Big].\label{FP2}
\end{align}
In particular, if $\alpha$ and $\beta$ are constants, we get
\begin{align}
    \frac{\partial^2}{\partial t \partial x} \mu_{t,x}= -\alpha D\mu_{t,x} + \tfrac{1}{2} \beta^2 D^2 \mu_{t,x}
    + tx \Big(\alpha^2 D^2 \mu_{t,x} +  \alpha \beta^2 D^3 \mu_{t,x} +\tfrac{1}{4}  \beta^4 D^4\mu_{t,x}\Big). \label{FP3}
\end{align}
\end{itemize}
\end{theorem}

\dproof The proof follows the idea of the proof of a similar result in \cite{AO}.\\
By the Dynkin formula we have, for a given locally bounded smooth function $\varphi$,
\begin{align*}
    &E[\varphi(Y(z))]-\varphi(y)=E\left[\int_{R_z} A\varphi(Y(\zeta)) d\zeta\right]
    \end{align*}
    where
    \begin{align}
    A\varphi(Y(\zeta))&= \varphi'(Y(\zeta)) \alpha(\zeta) + \tfrac{1}{2} \varphi''(Y(\zeta)) \beta^2(\zeta)\nonumber+ \int_{R_z} I(\zeta \bar{\wedge} \zeta^{\prime}) \Big\{\varphi''(Y(\zeta \vee \zeta'))\alpha(\zeta) \alpha(\zeta')\nonumber\\
    &+\tfrac{1}{2} \varphi^{(3)}(Y(\zeta \vee \zeta')) \big( \alpha(\zeta)\beta^2(\zeta')+\alpha(\zeta')\beta^2(\zeta)\big)\nonumber+\tfrac{1}{4} \varphi^{(4)}(Y(\zeta \vee \zeta')) \beta^2(\zeta) \beta^2(\zeta') \Big\} d\zeta'.
    \end{align}
Applying this to the function
$\varphi(y)=\varphi_w(y)=e^{-i wy}; \quad y,w\in \mathbb{R},$
where $i=\sqrt{-1}$,  and  with $\varphi'(y)=-iwe^{-i wy}$, $\varphi''(y)=-w^2e^{-i wy}$, $\varphi^{(3)}(y)=iw^3e^{-i wy}$, $\varphi^{(3)}(y)=w^4e^{-i wy}$  
we get
\begin{align}
    &E\left[ e^{-iwY(z)}\right]- e^{-iwy}=E\left[ \int_{R_z} \Big\{ -iw e^{-iwY(\zeta)} \alpha(\zeta) + \tfrac{1}{2}(-w^2)e^{-iwY(\zeta)} \beta^2(\zeta)\right.\nonumber\\
    &+\int_{R_z}I(\zeta \bar{\wedge} \zeta^{\prime}) \left(-w^2e^{-iwY(\zeta \vee \zeta')}\alpha(\zeta)\alpha(\zeta')\right.\nonumber+\tfrac{1}{2} i w^3 e^{-iwY(\zeta \vee \zeta')}\big(\alpha(\zeta) \beta^2(\zeta')+\alpha(\zeta')\beta^2(\zeta) \big) \nonumber\\
    &\left.\left.+\tfrac{1}{4}w^4 e^{-iwY(\zeta \vee \zeta')} \beta^2(\zeta) \beta^2(\zeta') \right ) d\zeta' \Big \} d\zeta\right].\label{2.8}
\end{align}
Note that
\begin{equation*}
     Y(\zeta)=Y(\zeta \vee \zeta')\mathbf{J}(\zeta,\zeta'), \quad 
    Y(\zeta')=Y(\zeta \vee \zeta')\mathbf{K}(\zeta,\zeta')
\end{equation*}

where
\begin{equation*}
    \mathbf{J}(\zeta,\zeta')= \mathbf{1}_{\zeta_1\geq \zeta_1' \& \zeta _2 \geq \zeta_2'} (\zeta,\zeta'), \quad  \mathbf{K}(\zeta,\zeta')= \mathbf{1}_{\zeta_1'\geq \zeta_1 \& \zeta' _2 \geq \zeta_2} (\zeta,\zeta').
\end{equation*}

Therefore \eqref{2.8} can be written
\begin{align}
   &E\left[ e^{-iwY(z)}\right]- e^{-iwy}=E\Big[ \int_{R_z} \int_{R_z} \Big\{ -iw e^{-iwY(\zeta \vee \zeta')} (\alpha(\zeta \vee \zeta')\mathbf{J} \nonumber\\
    &+ \tfrac{1}{2}(-w^2)e^{-iwY(\zeta \vee \zeta')} (\beta^2(\zeta \vee \zeta')\mathbf{J}\nonumber+\int_{R_z}I(\zeta \bar{\wedge} \zeta^{\prime}) \Big(-w^2e^{-iwY(\zeta \vee \zeta')}\alpha(\zeta \vee \zeta')\mathbf{J}\alpha(\zeta \vee \zeta')\mathbf{K}\nonumber\\
    &+\tfrac{1}{2} i w^3 e^{-iwY(\zeta \vee \zeta')}\big(\alpha(\zeta \vee \zeta')\mathbf{J} \beta^2(\zeta \vee \zeta')\mathbf{K}+\alpha(\zeta \vee \zeta')\mathbf{J}(\zeta,\zeta')\beta^2(\zeta \vee \zeta') \mathbf{K}(\zeta,\zeta')\big) \nonumber\\
    &\left.\left.+\tfrac{1}{4}w^4 e^{-iwY(\zeta \vee \zeta')} \beta^2(\zeta) \beta^2(\zeta') \right ) d\zeta' \Big \} d\zeta\right].\label{4.10}
\end{align}
Note that if $F$ denotes the Fourier transform operator, and $\mu_{\zeta}(dy)$ denotes the law of $Y(\zeta)$, we have
\begin{align*}
    E\left[ e^{-iwY(\zeta)}\right]=\int_{\R} e^{-iwy} \mu_{\zeta}(dy)=F[\mu_{\zeta}](w),
\end{align*}
and similarly
\begin{align*}
    E[e^{-iwY(\zeta \vee \zeta')}]=\int_{\R} e^{-iwy} \mu_{\zeta \vee \zeta'}(dy)=F[\mu_{\zeta \vee \zeta'}](w),
\end{align*}
Using this, and noting that $e^{-iwy}=F[\delta_y(\cdot)](w)$, we see that \eqref{2.8} can be written
\begin{align}
    F[\mu_z](w) - F[\delta_y(\cdot)](w)= &\int_{R_z} \Big\{ -iwF[\alpha \mu_{\zeta}](w) -\tfrac{1}{2} w^2 F[\beta^2 \mu_{\zeta}](w) \nonumber\\
 &+\int_{R_z} I(\zeta \bar{\wedge} \zeta^{\prime}) \Big( -w^2 F[\alpha(\zeta)\alpha(\zeta') \mu_{\zeta \vee \zeta'}](w) \nonumber\\
 & +\tfrac{1}{2}i w^3 F[\big(\alpha(\zeta)\beta^2(\zeta') +\alpha(\zeta')\beta^2(\zeta)\big)\mu_{\zeta \vee \zeta'}](w) \nonumber\\
 & +\tfrac{1}{4} w^4 F[\beta^2(\zeta) \beta^2(\zeta')\mu_{\zeta \vee \zeta'}](w) \Big) d\zeta' \Big\} d\zeta.\label{2.11}
\end{align}
Now we use that
\begin{equation*}
iwF[\alpha \mu_{\zeta}](w)= F[D(\alpha \mu_{\zeta})](w) \quad \text{and} \quad   -w^2 F[\beta^2 \mu_{\zeta}](w)=F[D^2({\beta^2 \mu_{\zeta})](w)}],  
\end{equation*}

where $D=\frac{\partial }{\partial y}, D^2=\frac{\partial ^2}{\partial y^2}$ denote derivatives with respect to $y$, in the sense of distribution.   Using (3.12) and (3.13), we can define terms with third and fourth-order derivatives.
\begin{equation*}
    -iw^3 F[\alpha\beta^2 \mu_{\zeta}](w)=F[D^3({\alpha\beta^2 \mu_{\zeta})](w)}] \quad \text{and} \quad  w^4 F[\beta^2\beta^2 \mu_{\zeta}](w)=F[D^4({\beta^2\beta^2 \mu_{\zeta})](w)}].
\end{equation*}

Then \eqref{2.11} can be written 
\begin{align*}
F[\mu_{\zeta} -\delta_y](w) &=F\Big[\int_{R_z}\Big \{-D[\alpha \mu_{\zeta}] + \tfrac{1}{2} D^2[\beta^2 \mu_{\zeta}]\nonumber\\
&+\int_{R_z}I(\zeta \bar{\wedge} \zeta^{\prime}) \Big( D^2[\alpha{(\zeta) \alpha(\zeta')} \mu_{\zeta \vee \zeta'}]-\tfrac{1}{2} D^3[\big(\alpha(\zeta)\beta^2(\zeta')+\alpha(\zeta') \beta^2 (\zeta)\big) \mu_{\zeta \vee \zeta'}]\nonumber\\
&+\tfrac{1}{4} D^4[\beta^2(\zeta) \beta^2(\zeta') \mu_{\zeta \vee \zeta'} \Big)d\zeta'\Big\}d\zeta \Big] (w).
\end{align*}
By uniqueness of the Fourier transform we conclude that $\mu$ satisfies the equation
\begin{align*}
    \mu_z = \delta_y + \int_{R_z}\Big\{A_1[\mu_{\zeta}] + \int_{R_z}L_{\zeta'}[\mu_{\zeta}] d\zeta' \Big \} d\zeta,
\end{align*}
which is \eqref{FP1}.\\
The rest of the theorem, i.e. \eqref{FP2} and \eqref{FP3} follow by Lemma \ref{4.1}.
\fproof

\begin{example} Consider the time -space SPDE
\small   
\begin{equation*}
    \frac{\partial^2}{\partial t \partial x} Y(t,x)=\beta \overset{\bullet}{B}(t,x); \quad t,x>0, \quad
Y(0,0)=y_0,
\end{equation*} 
where $\beta > 0$ is a constant.
By \eqref{FP1} the law $\mu_{t,x}(\cdot)=\mathcal{L}(Y(t,x))$ satisfies the integral equation
\small
\begin{align}
    \mu_z &= \delta_{y_0} + \int_{R_z}A_1[\mu_{\zeta}] d\zeta +\int_{R_z} \Big\{ \int_{R_z}L_{\zeta'}[\mu_{\zeta}] d\zeta' \Big \} d\zeta,\label{FPa}\nonumber\\
    &= \delta_{y_0} + \tfrac{1}{2} \beta^2\int_{R_z}  D^2[\mu_{\zeta}]d\zeta +\tfrac{1}{4} \beta^4 \int_{R_z} \Big\{\int_{R_z} I(\zeta \bar{\wedge} \zeta^{\prime}) D^4[\mu_{\zeta \vee \zeta'}]d\zeta' \Big\} d\zeta,
\end{align}
where as before we have put $z=(t,x)$.\\
As in Lemma \ref{4.1}, using the definition of $I(\zeta \bar{\wedge} \zeta^{\prime})$ and $\zeta \vee \zeta'$, we see that, with \newline $\zeta=(\zeta_1,\zeta_2), \zeta'=(\zeta'_1, \zeta'_2),$ we have
\small
\begin{align*}
    \int_{R_z} I(\zeta \bar{\wedge} \zeta^{\prime}) D^4[\mu_{\zeta \vee \zeta'}]d\zeta' =\zeta_2\int_{\zeta_1}^t D^4 [\mu_{(\zeta_1^{\prime},\zeta_2)}] d\zeta_1^{\prime}.
\end{align*}
In particular, if we assume that $\mu_{t,x}$ is absolutely continuous with respect to Lebesgue measure with Radon-Nikodym derivative $m_{(t,x)}(y)$, i.e., $\mu_{t,x}(dy)=m_{(t,x)}(y)dy$
then the equation \eqref{FP1} gets the following integro-differential form
\small
\begin{align}\label{FP2a}
m_{(t,x)}(y)&=m_{(0,0)}(y)+\tfrac{1}{2} \beta^2 \int_{R_z}  \frac{\partial ^2}{\partial y^2} m_{\zeta}(y)d\zeta \nonumber+\tfrac{1}{4} \beta^4 \int_{R_z} \Big\{\int_{R_z} I(\zeta \bar{\wedge} \zeta^{\prime}) \frac{\partial ^4}{\partial y^4}m_{\zeta \vee \zeta'}(y) d\zeta' \Big\} d\zeta\nonumber\\
&=\tfrac{1}{2} \beta^2\int_{R_z}  \frac{\partial ^2}{\partial y^2} m_{\zeta}(y)d\zeta +\tfrac{1}{4} \beta^4 \int_{R_z}\zeta_2\Big(\int_{\zeta_1}^t \frac{\partial ^4}{\partial y^4} m_{(\zeta_1^{\prime},\zeta_2)}(y) d\zeta_1^{\prime}\Big) d\zeta_1 d\zeta_2.
\end{align}
Differentiating this with respect to $t$ and $x$ we obtain\\
\small
\begin{align}
    \frac{\partial^2}{\partial x \partial t} m_{(t,x)}(y)&=\tfrac{1}{2}\frac{\partial^2}{\partial y^2} m_{(t,x)}(y)
    + \tfrac{1}{4}\frac{\partial}{\partial x}\left(\frac{\partial}{\partial t} \int_{R_z} \zeta_2\left(\int_{\zeta_1}^t \frac{\partial^4}{\partial y^4} m_{(\zeta_1',\zeta_2)}(y)d\zeta'\right)d\zeta_1 d\zeta_2\right)\nonumber\\
    &=\tfrac{1}{2}\frac{\partial^2}{\partial y^2} m_{(t,x)}(y)
    + \tfrac{1}{4}\frac{\partial}{\partial x}\left( \int_{R_z} \zeta_2 \frac{\partial^4}{\partial y^4} m_{(t,\zeta_2)}(y) d\zeta_1 d\zeta_2\right) \nonumber\\
    &=\tfrac{1}{2}\frac{\partial^2}{\partial y^2} m_{(t,x)}(y)
    + \tfrac{1}{4} \int_{0}^t x \frac{\partial^4}{\partial y^4} m_{(t,x)}(y) d\zeta_1  \nonumber\\
    &=\tfrac{1}{2}\frac{\partial^2}{\partial y^2} m_{(t,x)}(y)+ \tfrac{1}{4} t  x \frac{\partial^4}{\partial y^4} m_{(t,x)}(y). \label{fp2}
\end{align}
If we put $\beta=1$ the process $Y(t,x)$ is just the Brownian sheet, which is Gaussian and we know that in that case the law is absolutely continuous with respect to 2-dimensional Lebesgue measure with density $m_{(t,x)}(y)$ given by
\small
\begin{align*}
    m_{(t,x)}(y)=\frac{1}{\sqrt{2 \pi t x}} \exp\left\{ - \frac{(y-y_0)^2}{2tx}\right\}.
    \end{align*}
    \end{example}   
Below is the graph of the density of the Brownian sheet for a given $y_0=0.5.$

\begin{figure}[H]
\centering
\includegraphics[width=1\linewidth]{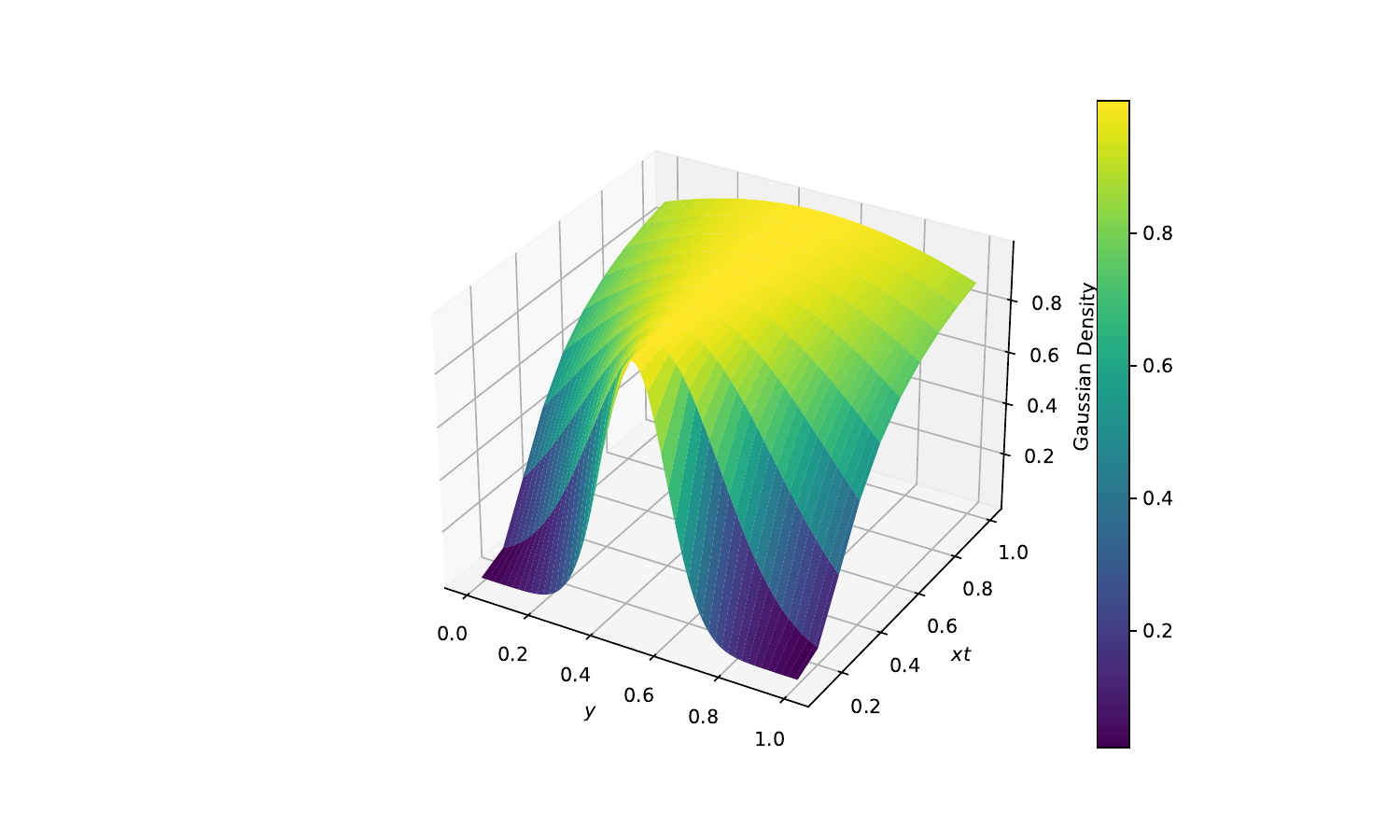}  
\caption{Brownian sheet density on a 2D grid }
\end{figure}

\begin{remark}
We now verify by explicit computation that this function $m_{(t,x)}(y)$ solves the Fokker-Planck equation \eqref{FP2}:\\
First we compute
\small
\begin{equation*}
     \frac{\partial}{\partial y}m_{(t,x)}(y)=- m_{(t,x)}(y) \frac{y-y_0}{tx}, 
\end{equation*}
and
\small
\begin{equation}\label{Y2}
         \frac{\partial ^2}{\partial y^2}m_{(t,x)}(y)
    =m_{(t,x)}(y)\left[-\frac{1}{tx}+\frac{(y-y_0)^2}{t^2 x^2}\right]
\end{equation}
and
\small
\begin{align*}
    \frac{\partial ^3}{\partial y^3}m_{(t,x)}(y)&=- m_{(t,x)}(y) \frac{y-y_0}{tx}\left[-\frac{1}{tx}+\frac{(y-y_0)^2}{t^2 x^2}\right]+m_{(t,x)}(y) \frac{2(y-y_0)}{t^2x^2}\nonumber\\
&=m_{(t,x)}(y)\left[-\frac{(y-y_0)^3}{t^3x^3}+\frac{3(y-y_0)}{t^2 x^2}\right],
\end{align*}
and
\small
\begin{align}
    \frac{\partial ^4}{\partial y^4}m_{(t,x)}(y)&=- m_{(t,x)}(y) \frac{y-y_0}{tx}\left[-\frac{(y-y_0)^3}{t^3x^3}+\frac{3(y-y_0)}{t^2 x^2}\right]\nonumber\\
&+m_{(t,x)}(y) \left[\frac{-3(y-y_0)^2}{t^3x^3}+\frac{3}{t^2 x^2}\right]\nonumber\\
%&=m_{(t,x)}(y)[\frac{(y-y_0)^4}{t^4x^4}-\frac{6(y-y_0)^2}{t^3 x^3}+\frac{3}{t^2 x^2}]\nonumber\\
&=\frac{1}{t^2 x^2} m_{(t,x)}(y)\left[\frac{(y-y_0)^4}{t^2x^2}-\frac{6(y-y_0)^2}{t x}+3\right].\label{Y4}
\end{align}
Next, we compute
\small
\begin{align*}
    \frac{\partial }{\partial t}m_{(t,x)}(y)&=m_{(t,x)}(y)\left[-\frac{1}{2t} +\frac{(y-y_0)^2}{2xt^2}\right],
\end{align*}
and
\small
\begin{align}
    \frac{\partial ^2}{\partial t \partial x}m_{(t,x)}(y)&=m_{(t,x)}(y)\left[-\frac{1}{2x} +\frac{(y-y_0)^2}{2x^2t}\right]\left[-\frac{1}{2t} +\frac{(y-y_0)^2}{2xt^2}\right]\nonumber\\
    &+m_{(t,x)}(y)\left(-\frac{(y-y_0)^2}{2x^2t^2}\right)\nonumber\\
    %&=m_{(t,x)}(y)[\frac{1}{4tx}-2\frac{(y-y_0)^2}{4x^2t^2}+\frac{(y-y_0)^4}{4x^3t^3}-\frac{(y-y_0)^2}{2x^2t^2}]\nonumber\\
    &=m_{(t,x)}(y)\left[\frac{1}{4tx}-\frac{(y-y_0)^2}{x^2t^2}+\frac{(y-y_0)^4}{4x^3t^3}\right].
    \label{fp1}
\end{align}
By substituting \eqref{Y2} and \eqref{Y4} into \eqref{fp1} we get \eqref{fp2}, as required.

\end{remark}

%%%%%%%%%%%%%%%%%%%%%%%%%%%%%%%%%%%%%%%%%%%%%%%%%%%%%%%%%%%%%%%%%%%%%%%%%%%%%%%%%%%%%%%%%%%%%%%%%%%%%%%%%%%%%%%%%%%%%%%%%%%%%%%%%%%%%%%%%%%%%%%%%%%%%%%%%%%%%%%%%%%%%

\end{document}